\baselineskip=14pt
\parskip=10pt

\def\epsilon{\varepsilon}

\font\eightrm=cmr8  
\font\eighttt=cmtt8
\magnification=\magstephalf

\parindent=0pt
\overfullrule=0in
 
\bf
\noindent
THEOREMS FOR A PRICE: Tomorrow's Semi-Rigorous Mathematical Culture 
\medskip
\noindent
\it
\hskip150pt \relax Doron Zeilberger\footnote{$^1$}
{\eightrm 
Department of Mathematics, Temple University,
Philadelphia, PA 19122, USA. 
{\eighttt zeilberg@math.temple.edu .}
Supported in part by the NSF. Based in part on a Colloquium talk given
at Rutgers University. Appeared in Notices of the Amer. Math. Soc.
40(1993), 978-981. Reprinted in Math. Intell. 16(4) (1994), 11-14.
}
\bigskip
\rm
{ \bf Today}
 
The most fundamental precept of the mathematical faith is
{\it thou shalt prove everything rigorously.}
While the practitioners of {\it mathematics} differ in their
view of what constitutes a rigorous proof, and there are
fundamentalists who insist on even a more rigorous rigor
than the one practiced by the mainstream, the belief
in this principle could be taken as the {\it defining property} of
{\it mathematician}.
 
{\bf The Day After Tomorrow}
 
There are writings on the wall that, now that the 
silicon savior has arrived, a new testament is going to be written. Although
there will always be a small
group of ``rigorous'' old-style mathematicians(e.g. [JQ]) who will insist that
the true religion is theirs,
and that the computer is a false Messiah,
they may be viewed by future mainstream mathematicians as 
a fringe sect of harmless eccentrics, like mathematical
physicists are viewed by regular physicists today.
 
The computer has already started doing to mathematics what the
telescope and microscope did to astronomy and biology.
In the future, not all mathematicians will care about absolute certainty, 
since there will be so many exciting new facts to discover: mathematical
pulsars and quasars that  will make the Mandelbrot set seem like a mere 
Jovian moon. We will have (both human and machine\footnote{$^2$}
{\eightrm For example, my computer Shalosh B. Ekhad,
and its friend Sol Tre, already have a
non-trivial publication list, e.g. [E], [ET]}
) professional {\it theoretical} mathematicians, who will
develop conceptual paradigms to make sense out of the empirical data,
and who will reap Fields medals along with (human and machine)
{\it experimental} mathematicians. Will there  still be a place for
{\it mathematical} mathematicians?
 
This will happen after a transitory age of {\it semi-rigorous mathematics},
in which identities (and perhaps other kinds of theorems) will carry
price-tags.
 
{\bf A Taste Of Things To Come}
 
To get a glimpse of how mathematics will be practiced in the not too
far future, I will describe the case of algorithmic proof theory  for
{\it hypergeometric identities}[WZ*][Z*][Ca]. In this theory, it is possible to
rigorously prove, or refute, any conjectured identity belonging to a wide class
of identities, that includes most of the identities between the classical
special functions of mathematical physics.
 
Any such identity is proved by exhibiting a {\it proof certificate},
that reduces the proof of the given identity to that of a finite
identity among rational functions, and hence, by clearing denominators,
to that between specific polynomials.
 
This algorithm can be performed successfully on all
``natural identities'' we are now aware of.
It is easy, however, to concoct artificial examples
for which the running time, and memory, are prohibitive. Undoubtedly,
in the future ``natural'' identities will be encountered whose complete
proof will turn out to be not worth the money. We will see later, how,
in such cases, one can get ``almost certainty'' with a tiny fraction
of the price, along with the assurance that if we robbed a bank, we would be
able to know for sure. 
 
This is vaguely reminiscent of {\it transparent proofs} introduced
recently in theoretical computer science[Ci][ALMSS][AS].
The result that there exist short theorems
having arbitrarily long proofs, a consequence of G\"odel's incompleteness
theorem, also comes to mind [S].\footnote{$^3$}
{\eightrm Namely, the ratio (proof length)/(theorem length) grows fast
enough to be non-recursive. 
Adding an axiom can shorten proofs by recursive amounts [G], [D].}
I speculate that similar developments
will occur elsewhere in mathematics, and will ``trivialize'' large
parts of mathematics, by reducing mathematical truths to routine,
albeit possibly very long, and exorbitantly expensive to check,
{\it ``proof certificates''}. These proof certificates would also 
enable us, by plugging in random values,
to assert ``probable truth'' very cheaply.
 
{\bf Identities}
 
Many mathematical theorems are {\it identities}: statements of
type ``='', which take the form $A=B$. Here is a sampler, in roughly
an increasing order of sophistication.
 
{\bf 1.} $2+2=4$.
 
{\bf 2.} $(a+b)^3=a^3+3 a^2 b + 3a b^2 + b^3$.
 
{\bf 3.} $sin(x+y)=sin(x)cos(y)+cos(x)sin(y)$.
 
{\bf 4.} $F_{n+1} F_{n-1}-F_{n}^2=(-1)^{n+1} \,$.
 
{\bf 5.} $(a+b)^n = \sum_{k=0}^{n} {n \choose k} a^k b^{n-k} \,$.
 
{\bf 6.} $\sum_{k=-n}^{n} (-1)^k {{2n} \choose {n+k} }^3 = {{3n} \choose {n}}$.
 
{\bf 7.} Let $(q)_r := (1-q)(1- q^2 ) \dots (1-q^r)$, then
 
$$
\sum_{r=0}^{n} 
{
{ q^{r^2} } \over {(q)_r (q)_{n-r} }
}
=
\sum_{r=-n}^{n}
{
{(-1)^r q^{(5 r^2 -r)/2}}
\over
{ (q)_{n-r} (q)_{n+r} }
} \quad .
$$
 
{\bf 7'.} Let $(q)_r$ be as in 7.
 
$$
\sum_{r=0}^{\infty} 
{
{ q^{r^2} } \over {(q)_r }
}
=
\prod_{i=0}^{\infty} (1-q^{5i+1})^{-1} (1-q^{5i+4})^{-1} \quad .
$$
 
{\bf 8.} Let $H_n$ be given by,
 
$$
H_n = H_n (q) =
{ {(1+q)} \over {(1-q)} }
{ {(1+q^2)} \over {(1-q^2)} }
\dots
{ {(1+q^n)} \over {(1-q^n)} } \quad ,
$$
 
then
 
$$
\left ( {\sum_{k=0}^{n} {{2 (- q^{n+1} )^k } \over {1+ q^k}}
{H_k} } \right )^4 \,
\sum_{k=-n}^n 
{{4 (-q)^k } \over {(1+ q^k )^2 }} { { H_{n+k} } \over { H_n} }
{{ H_{n-k} } \over { H_n }}
=  ( \sum_{k= -n }^n (-q)^{ k^2 } )^4 \; .
$$
 
{\bf 8'.}
 
$$
( \sum_{k= - \infty }^{\infty} q^{ k^2 } \, )^4 =
1+ 8  \sum_{k=1}^{\infty} {{q^k} \over { (1 + (-q)^k )^2 }} \; .
$$
 
{\bf 9.} Analytic Index=Topological Index.
 
{\bf 10.} $Re(s)= {{1} \over {2}} $ for every 
non-real $s$ such that $\zeta (s) =0$.
 
All the above identities are trivial, except possibly the last
two, which I think quite likely will be considered trivial in
two hundred years. I will now explain.
 
{\bf Why are the first 8 identities trivial?}
 
The first identity, while {\it trivial} nowadays, was very deep
when it was first discovered, independently, by several anonymous
cave-dwellers.
It is a general, abstract theorem, that contains, as special cases,
many apparently unrelated theorems: Two bears and Two bears make
Four bears, Two apples and Two apples make Four apples, etc.
It was also realized that in order to prove it rigorously, it suffices
to prove it for any one special case, say, marks on the
cave's wall.
 
The second identity: $(a+b)^3=a^3+3 a^2 b + 3a b^2 + b^3$,  is of one
level of generality higher. Taken literally (in the {\it semantic} sense
of the word {\it literally}), it is a fact about {\it numbers}. For any
specialization of $a$ and $b$ we get yet another correct numerical
fact, and as such it requires a ``proof'', invoking the commutative,
distributive and associative ``laws''. However, 
it is completely routine when viewed {\it literally},
in the syntactic sense, i.e. in which $a$ and $b$ are no longer
symbols {\it denoting} numbers, but rather represent themselves,
qua (commuting) literals. 
This shift in emphasis roughly corresponds to the transition from
Fortran to Maple, i.e. from {\it numeric} computation to
{\it symbolic} computation.
 
Identities {\bf 3} and {\bf 4} can be easily embedded in classes of routinely
verifiable identities in several ways. One way is by
{\it defining} $cos(x)$ and $sin(x)$ by $(e^{ix} + e^{-ix})/2$
and $(e^{ix} - e^{-ix})/(2i)$, and the
Fibonacci numbers $F_n$ by Binet's formula.
 
Identities {\bf 5-8} were, until recently, considered genuine non-trivial
identities, requiring a human demonstration. One particularly nice
human proof of {\bf 6}, was given by Cartier and Foata[CF].
A one-line computer-generated proof of identity {\bf 6} is given in
[E]. Identities {\bf 7} and {\bf 8}
are examples of so-called {\it q-binomial coefficient
identities}, (a.k.a. {\it terminating q-hypergeometric series}.)
All such identities are now routinely provable[WZ2](see below.)
The machine-generated proofs of {\bf 7} and {\bf 8} appear in
[ET] and [AEZ] respectively. Identities
{\bf 7} and {\bf 8} immediately imply, 
by taking the limit $n \rightarrow \infty$,
identities {\bf 7'} and {\bf 8'}, which, in turn, are equivalent to two famous
number-theoretic statements: The first Rogers-Ramanujan identity,
that asserts that the number of partitions of an integer into parts
that leave remainder $1$ or $4$ when divided by $5$ equals
the number of partitions of that integer into parts that differ
from each other by at least $2$, and Jacobi's
theorem that asserts that
the number of representations of an integer as a sum of
$4$ squares, equals $8$ times the sum of its divisors that
are not multiples  of $4$.
 
{\bf The WZ Proof Theory}
 
Identities {\bf 5-8} involve sums of the form
 
$$
\sum_{k=0}^{n} F(n,k)  \quad ,
\eqno(Sum)
$$
 
where the summand, $F(n,k)$, is a
{\it hypergeometric term}, (in {\bf 5}, and {\bf 6}), or a
{\it q-hypergeometric term}, (in {\bf 7}, and {\bf 8}),
in both $n$ and $k$, which means
that both quotients $F(n+1,k)/F(n,k)$ and $F(n,k+1)/F(n,k)$ are
{\it rational functions} of $(n,k)$ ($(q^n,q^k,q)$ respectively).
 
For such sums, and multi-sums, we have ([WZ2]) the following result.
 
{\bf The Fundamental Theorem of Algorithmic Hypergeometric Proof Theory:}
 
Let $F(n; k_1 , \dots , k_r)$ be a {\it proper} (see [WZ2])
hypergeometric term in 
all of $(n; k_1 , \dots , k_r )$.
Then there exist polynomials
$p_0(n), \dots , p_L(n)$ and {rational functions}  
$R_j(n; k_1 , \dots, k_r )$ such that $G_j:=R_j F$ satisfy
 
$$
\sum_{i=0}^{L}
p_i(n)F(n+i; k_1 , \dots , k_r )=
$$
$$
 \sum_{j=1}^r [G_j(n; k_1 , \dots , k_j +1 , \dots , k_r)-
G_j(n; k_1 , \dots , k_j , \dots , k_r)]
\eqno (multiWZ)
$$
 
and hence, if for every specific $n$, $F(n; -)$ has compact support in
$(k_1, \dots , k_r)$, the definite sum, $g(n)$, given by
 
$$
g(n):=\sum_{k_1, \dots , k_r} F(n; k_1 , \dots , k_r )
\eqno(multiSum)
$$
 
satisfies the linear recurrence equation with polynomial coefficients:
 
$$
\sum_{i=0}^{L} p_i(n) g(n+i) =0 \quad .
\eqno(P-recursive)
$$
 
$(P-recursive)$ follows from $(multiWZ)$ by summing over 
$\{ k_1 , \dots , k_r \}$, and observing that all the sums
on the right telescope to zero.
 
If the recurrence happens to be first order, i.e. $L=1$ above,
then it can be written in {\it closed form}: for example the solution
of the recurrence $(n+1)g(n)-g(n+1)=0$, $g(0)=1$, is $g(n)=n!$.
 
This ``existence'' theorem also implies an algorithm for finding
the recurrence (i.e. the $p_i$) and the accompanying certificates
$R_j$, see below.
 
An analogous theorem holds for q-hypergeometric series[WZ2][K].
 
Since we know how to find, and prove, the recurrence satisfied
by any given hypergeometric sum or multi-sum, we have an
effective way of proving any equality of two such sums, or
the equality of a sum with a conjectured sequence. All we
have to do is check whether both sides are solutions of the
same recurrence, and match the appropriate number of initial values.
Furthermore, we can also use the algorithm to find new identities.
If a given sum yields a first-order recurrence,
it can be solved, as mentioned above, and the sum in question
turns out to be explicitly evaluable. If the recurrence obtained
is of higher order, then most likely the sum is not explicitly-evaluable
(in closed form), and Petkovsek's algorithm[P], that
decides whether a given linear recurrence (with polynomial
coefficients) has {\it closed form}
solutions, can be used to find out for sure.
 
{\bf Almost Certainty For An $\bf{\epsilon}$ Of The Cost }
 
Consider identity $(multiSum)$ once again, where $g(n)$ is ``nice''.
Dividing through by $g(n)$ and letting $F \rightarrow F/g$ (Herb Wilf's
wonderful trick), we can assume that we have to prove an identity of
the form
 
$$
\sum_{k_1, \dots , k_r} F(n; k_1 , \dots , k_r ) \, = \, 1 \quad .
\eqno(Nice)
$$
 
The WZ theory promises you that the left side satisfies some
linear recurrence, and if the identity is indeed true, then the
sequence $g(n)=1$ should be a solution (in other words
$p_0 (n) + \dots + p_L (n) \equiv 1$). For the sake of
simplicity, let's assume that the recurrence is minimal,
i.e. is $g(n+1)-g(n)=0$,
(this is true anyway in the vast majority of the cases.)
To prove the identity, by this method, we have to find {\it rational functions}
$R_j(n;k_1 , \dots , k_r)$ such that the $G_j := R_j F$ satisfy:
 
$$
F(n+1; k_1 , \dots , k_r ) - F(n; k_1 , \dots , k_r )=
$$
$$
 \sum_{j=1}^r [G_j(n; k_1 , \dots , k_j +1 , \dots , k_r)-
G_j(n; k_1 , \dots , k_j , \dots , k_r)] \, .
\eqno (multiWZ')
$$
 
By dividing $(multiWZ')$ through by $F$, and clearing denominators,
we get a certain functional equation for the $R_1 , \dots , R_r$,
from which it is possible to determine their denominators 
$Q_1 , \dots , Q_r$. Writing $R_j=P_j/Q_j$, the proof boils down
to finding {\it polynomials} $P_j ( k_1 , \dots , k_r )$, with
coefficients that are rational functions in $n$ and possibly
other (auxiliary) parameters. It is easy to predict upper
bounds for the degrees of the $P_j$ in $(k_1 , \dots , k_r)$.
We then express each $P_j$ symbolically, with ``undetermined''
coefficients, and substitute into the above-mentioned functional equation.
We then expand,
and equate coefficients of all monomials $k_1^{a_1} \dots k_r^{a_r}$,
and get an (often huge) system of inhomogeneous linear
equations with {\it symbolic} coefficients. The proof boils
down to proving that this inhomogeneous system of linear
equations has a solution. It is very time-consuming to solve
a system of linear equations with {\it symbolic} coefficients.
By plugging in
specific values for $n$ and the other parameters, if present,
one gets a system with {\it numerical} coefficients, which is much
faster to handle. Since it is unlikely that a random system of 
inhomogeneous linear
equations with more equations than unknowns can be solved,
the solvability of the system for a number of special values of
$n$ and the other parameters
is a very good indication that the identity is indeed true.
It is a waste of money to get absolute certainty,
unless the conjectured identity in question is known
to imply the Riemann Hypothesis.
 
{\bf Semi-Rigorous Mathematics  }
 
As wider classes of identities, and perhaps even other kinds of
classes of theorems, become routinely-provable, we might witness
many results for which we would know how to find a proof (or refutation),
but we would be unable, or unwilling, to pay for finding such proofs,
since ``almost certainty'' can be bought so much cheaper. I can envision
an abstract of a paper, c. 2100, that reads :
``We show, in a certain precise sense,
that the Goldbach conjecture is true with probability larger than $0.99999$,
and that its complete truth could be determined with a budget of 
$\$ 10B$.''
 
It would be then OK to rely on such a priced theorem, provided that
the price is stated explicitly. Whenever statement $A$, whose price is
$p$, and statement $B$, whose price is $q$, are used to deduce statement $C$, 
the latter becomes a priced theorem priced at $p+q$.
 
If a whole chain of boring identities would turn out to imply
an interesting one, we might be tempted to redeem all these intermediate
identities, but we would not be able to buy out the whole store,
and most identities would have to stay unclaimed.
 
As absolute truth becomes more and more expensive, we would sooner or later
come to grips with the fact that few non-trivial results could be
known with old-fashioned certainty. 
Most likely we will wind up abandoning
the task of keeping track of price altogether, and complete the metamorphosis
to non-rigorous mathematics.
 
{\bf Note:} Maple programs for proving hypergeometric identities
are available by anonymous {\tt ftp} to {\tt math.temple.edu}, in
directory {\tt pub/zeilberger/programs}. a Mathematica implementation
of the single-summation program,
can be obtained from Peter Paule at {\tt paule@risc.uni-linz.ac.at}\quad .
 
\medskip
{\bf References}
 
[AZ] G. Almkvist and D. Zeilberger,
{\it The method of differentiating under the
integral sign}, J. Symbolic Computation {\bf 10}(1990), 571-591.
 
[AEZ] G. E. Andrews, S. B. Ekhad, and D. Zeilberger {\it A short proof of
Jacobi's formula for the number of representations of an integer
as a sum of four squares}, Amer. Math. Monthly, {\bf 100}(1993), 274-276.
 
[ALMSS] S. Arora, C. Lund, R. Motwani, M. Sudan, and M. Szegedy,
{\it Proof verification and intractability of approximation problems},
Proc. 33rd Symp. on Foundations of Computer Science (FOCS), IEEE Computer
Science Press, Los Alamitos, 1992, pp. 14-23.
 
[AS] S. Arora, and M. Safra, {\it Probabilstic checking of proofs},
ibid, pp. 2-13.
 
[Ca] P. Cartier, {\it D\'emonstration ``automatique'' d'identit\'es et
fonctions hyperg\'eometriques [d'apres D. Zeilberger]}, S\'eminaire
Bourbaki, expos\'e $n^{o}$ $746$, Ast\'erisque {\bf 206}, $41-91$, SMF, 1992.
 
[CF] P. Cartier and D. Foata, {\it Probl\`emes combinatoires de commutation
et r\'earrangements}, Lecture Notes Math. v. {\bf 85}, Springer 1969.
 
[Ci] B. Cipra, {\it Theoretical computer scientists develop transparent proof
techniques}, SIAM News, {\bf 25/3}(May 1992).
 
[D] J. Dawson, {\it The G\"odel incompleteness theorem from a length of
proof perspective}, Amer. Math. Monthly {\bf 86}(1979), 740-747.
 
[E] S. B. Ekhad, {\it A very short proof of Dixon's theorem}, J. Comb. Theo.
Ser. A {\bf 54}(1990), 141-142.
 
[ET] S. B. Ekhad and S. Tre, {\it A purely verification proof of the
first Rogers-Ramanujan identity}, J. Comb. The. Ser. A {\bf 54}(1990),
309-311.
 
[G] K. G\"odel, {\it On length of proofs},
(in German) {\it Ergeb. Math. Colloq.} {\bf 7}(1936), 23-24. Translated
in ``The Undecidable'' (M. Davis, Ed.), Raven Press, Hewitt, NY, 1965, 82-83.
 
[JQ] A. Jaffe and F. Quinn, {\it ``Theoretical mathematics'': Toward
a cultural synthesis of mathematics and theoretical physics}, Bulletin(N.S.)
of the Amer. Math. Soc. {\bf 29}(1993), 1-13.
 
[K] T. H. Koornwinder, {\it Zeilberger's algorithm and its q-analogue},
University of Amsterdam Mathematics department preprint.
 
[P] M. Petkovsek, {\it Hypergeometric solutions of linear recurrence
equations with polynomial coefficients}, J. Symbolic Computation
{\bf 14}(1992), 243-264.
 
[S] J. Spencer, {\it Short theorems with long proofs}, Amer. Math. Monthly
{\bf 90}(1983), 365-366.
 
[WZ1] H. S. Wilf and D. Zeilberger,
{\it Rational functions certify combinatorial
identities}, J. Amer. Math. Soc. {\bf 3}(1990), 147-158.
 
[WZ1a] ------,{\it Towards computerized proofs of identities},
Bulletin(N.S.) of the Amer. Math. Soc. {\bf 23}(1990), 77-83.
 
[WZ2] -------, {\it An algorithmic proof theory for hypergeometric
(ordinary and ``q'') multisum/integral identities}, Invent. Math. 
{\bf 108}(1992), 575-633.
 
[WZ2a] ------, {\it Rational function certification of
hypergeometric multi-integral/sum/``q'' identities}, 
Bulletin(N.S.) of the Amer. Math. Soc. {\bf 27}(1992) 148-153.
 
[Z1] D. Zeilberger, {\it A Holonomic systems approach to special functions
identities}, J. of Computational and Applied Math. {\bf 32}(1990),
321-368.
 
[Z2]------, {\it  A Fast Algorithm for proving terminating hypergeometric
identities}, Discrete Math {\bf 80}(1990), 207-211.
 
[Z3]------, {\it The method of creative telescoping}, J. Symbolic Computation
{\bf 11}(1991), 195-204.
 
[Z4]------, {\it Closed Form (pun intended!)}, 
``Special volume in memory of Emil Grosswald'', M. Knopp
and M. Sheingorn,
eds., Contemporary Mathematics {\bf 143}, 579-607, AMS, Providence, 1993.
 
\medskip
 
July 1993 ; Revised: August 1993.
 
\bye